\documentclass{amsart}

\usepackage{amssymb}
\usepackage{color}
\usepackage{amsxtra} 
\usepackage{mathrsfs} 
\usepackage{yfonts} 
\usepackage[dvipdfmx]{hyperref}

\numberwithin{equation}{section}

\newtheorem{theorem}{Theorem}[section]

\theoremstyle{definition}

\theoremstyle{remark}
\newtheorem{remark}{Remark}

   \makeatletter
\renewcommand*{\eqref}[1]{%
  \hyperref[{#1}]{\textup{\tagform@{\ref*{#1}}}}%
}
\makeatother

\allowdisplaybreaks[3]

\begin{document}

\title[inviscid shear flow]
{Remark on a lower bound of  perturbations in 2D inviscid shear flow in a periodic box} 

\author{Tsuyoshi Yoneda} 
\address{Graduate School of Mathematical Sciences, University of Tokyo, Komaba 3-8-1 Meguro, 
Tokyo 153-8914, Japan } 
\email{yoneda@ms.u-tokyo.ac.jp} 


\date{\today} 



\begin{abstract} 
Bedrossian and Masmoudi (2015) showed nonlinear stability of 
2D inviscid parallel shear flows in the infinite stripe $(\mathbb{R}/2\mathbb{Z})\times\mathbb{R}$.
Inspired by their result, combining a diffeomorphism group result by Misio{\l}ek (1993),
 in this paper we give a lower bound of  perturbations in 2D inviscid shear flow in a periodic box.
\end{abstract} 

\maketitle

\section{Introduction and the main theorem} 
\label{sec:Intro}

For $m\in\mathbb{N}$,
let  us define a stripe type of periodic box as $\mathbb{T}^2:=(\mathbb{R}/2\mathbb{Z})\times(\mathbb{R}/(2m\mathbb{Z}))$.
The incompressible Euler equations in $\mathbb{T}^2$ are described as follows:
\begin{eqnarray}\label{Euler}
&  &\partial_t u +(u \cdot \nabla) u
  =- \nabla
p,\ x=(x_1,x_2)\in \mathbb T^2,\ t>0, \\
\nonumber
& & \nabla \cdot u=0, \ u|_{t=0}=u_0,
\end{eqnarray}
where  $u=u(t,x)$ is the velocity and $p=p(t,x)$ is the pressure.
Throughout this paper we always handle smooth solutions, that is, we take initial velocity
$u_0\in C^\infty(\mathbb{T}^2)$.
It is well known that the corresponding unique smooth solution $u$ exists globally in time (see \cite{MB} for example).
Next we define the Lagrangian flow $\eta$. For $u\in C^\infty([0,\infty)\times \mathbb{T}^2)$, let $\eta(t,x)$ be a solution to the following ODE:
\begin{equation*}
\frac{d}{dt}\eta(t,x) = u(t, \eta(t,x))=:u\circ \eta,\quad 
\eta(0,x) =:e(x)=x. 
\end{equation*}
We can express it in Riemannian geometric language using the exponential map of the right-invariant metric (defined by the $L^2$ kinetic energy of the fluid) on the group $\mathscr{D}^s_\mu(\mathbb{T}^2)$ of volume-preserving diffeomorphisms of $\mathbb{T}^2$ of Sobolev class $H^s$ with $s>2$, namely $\eta(t) = \exp_e tu_0$ (see for example \cite{EM} or \cite{M}).
Clearly, $\eta$ uniquely exists globally in time and is smooth.
Now let us define the inviscid parallel shear flow (stationary Euler flow): 
\begin{equation*}
u_\infty(x)=(f_m(x_2), 0), \quad f_m\in C^\infty(\mathbb{R}/(2m\mathbb{Z})).
\end{equation*}
Here we choose $f_m$ in the following way:
First let us take $f\in C^\infty(\mathbb{R})$ with $\|f'\|_{L^\infty}<\infty$,
and then 
we take $f_m$ satisfying
\begin{equation}\label{f_m}
\begin{cases}
&\displaystyle\lim_{m\to\infty}\sup_{-\frac{m}{2}<x_2\leq\frac{m}{2}}\left(|f_m(x_2)-f(x_2)|+|f'_m(x_2)-f'(x_2)|\right)=0,\\
&\displaystyle\sup_m\|f'_m\|_{L^\infty}<\infty.
\end{cases}
\end{equation}

For any $U_{in}\in C^\infty$ with $\nabla\cdot U_{in}=0$, 
  from now on, 
we consider the following family of $L^2$ geodesics in $\mathscr{D}^s_\mu$, or equivalently, fluid flows in $\mathbb{T}^2$, starting at the identity configuration $e$ in the direction $u_\infty + \sigma U_{in}$:
\begin{equation*}
\tilde\eta(\sigma,t):=\exp_e(t(u_\infty+\sigma U_{in})),\quad \sigma\in(-1,1),\quad t\geq 0.
\end{equation*}
 Note that $\partial_\sigma\tilde\eta(\sigma,t)|_{t=0}=0$.
Let $\tilde u=\tilde u(\sigma,t)$ be the corresponding solution to the Euler equations \eqref{Euler} for the initial data $u_\infty+\sigma U_{in}$, and let us set
\begin{equation*}
\text{perturbation:}\quad\quad U=(U_1,U_2)= U(\sigma,t):=\tilde u(\sigma,t)-u_\infty.
\end{equation*}
 Note that $U(0,t)=0$ $(t>0)$.
In the infinite stripe case $(\mathbb{R}/2\mathbb{Z})\times \mathbb{R}$ (formally  $m=\infty$),
 Bedrossian and  Masmoudi \cite{BM} showed the following nonlinear stability result (we rewrite their statement into our setting): 
 For any $U_{in}\in \mathcal G^{\lambda_0}$ with $\|U_{in}\|_{\mathcal G^{\lambda_0}}\sim 1$ (for the precise definition of the Gevrey-norm $\|\cdot\|_{\mathcal G^{\lambda_0}}$, see \cite{BM}. In this paper, we just use the embeddings: $\|\cdot\|_{L^2}\leq \|\cdot\|_{\mathcal G^{\lambda}}$ and $\|\nabla \cdot\|_{L^\infty}\leq \|\cdot\|_{H^3}\leq \|\cdot\|_{\mathcal G^{\lambda}}$) and $\sigma>0$, then there is $\tilde f$ sufficiently smooth and 
appropriately decaying (depending on $\sigma$, but controlled by $\|\tilde f'\|_{L^\infty}\leq \|\tilde f\|_{\mathcal G^{\lambda'}}<\sigma\ll 1$), such that if  
\begin{equation*}
u_\infty(x)=
\begin{pmatrix}
f(x_2)\\
0
\end{pmatrix}
\quad\text{with}\quad f(x_2)=x_2+\tilde f(x_2),
\end{equation*}
 then we have the following stability estimate:
\begin{equation}\label{BM-estimate}
\|U_1(\sigma,t)
\|_{L^2}\lesssim\frac{\sigma}{\langle t\rangle} 
\quad\text{and}\quad
\|U_2(\sigma,t)\|_{L^2}\lesssim\frac{\sigma}{\langle t\rangle^2},
\end{equation} 
where $\langle t\rangle:=(1+|t|^2)^{1/2}$.

On the other hand, inspired by their result, we constructed the following lower bound.
\begin{theorem}
Let us consider the perturbation $U$  of the shear flow $(f_m,0)$ in $(\mathbb{R}/2\mathbb{Z})\times(\mathbb{R}/(2m\mathbb{Z}))$ satisfying \eqref{f_m}.
For any $\epsilon>0$ and 
 $T>0$, there is  $\sigma_0>0$ such that  if $\sigma<\sigma_0$, we have 
\begin{equation*}
\|U_1(\sigma,t)\|_{L^2(\mathbb{T}^2)}\gtrsim \sigma\quad\text{or}\quad \|U_2(\sigma,t)\|_{L^2(\mathbb{T}^2)}\gtrsim \frac{\sigma}{\langle t\rangle^{1+\epsilon}}
\end{equation*}
for  $m\in\mathbb{N}$ and $t\in[0,T]$. In particular, this lower bound is independent of $m$ and $T$.
\end{theorem}

\begin{remark}
If the energy of the initial velocity $u_0=u_\infty+\sigma U_{in}$
is not equal to the energy of the shear flow, namely, $\|u_0\|_{L^2}\not=\|u_\infty\|_{L^2}$,
then we can easily see that the energy of
perturbation $\|U(\sigma,t)\|_{L^2}$ is not zero for any $t>0$. 
This is just due to the triangle inequality (the following is the case when $\|u_0\|_{L^2}>\|u_\infty\|_{L^2}$):
\begin{equation*}
\|u_0\|_{L^2}=\|\tilde u(\sigma,t)\|_{L^2}\leq \|u_\infty\|_{L^2}+\|U(\sigma,t)\|_{L^2}
\Rightarrow \|u_0\|_{L^2}-\|u_\infty\|_{L^2}\leq \|U(\sigma,t)\|_{L^2}.
\end{equation*}
However, this estimate cannot exclude the case $\|u_0\|_{L^2}=\|u_\infty\|_{L^2}$ (with $u_0\not=u_\infty$).
\end{remark}

\section{Proof of the main theorem}

First we show $\partial_\sigma^k U(\sigma,\cdot)\in C^\infty([0,\infty)\times\mathbb{T}^2)$ ($\sigma\in(-1,1)$, $k=1,2,\cdots$).
In order to do so, just multiply $\partial_\sigma^k$ to \eqref{Euler} on both sides (replace $u$ with $\tilde u$), and just apply inductively the existence theory of the usual Euler equations (for the existence theory itself, see \cite{MB} example).
Since
$\tilde u(\sigma,t)=u_\infty+ U(\sigma, t)$, we see
\begin{equation*}
\partial_t\tilde\eta=u_\infty\circ\tilde\eta+ U\circ\tilde\eta.
\end{equation*}
Thus by differentiation in $\sigma$ on both sides, we have 
\begin{equation*}
\partial_t\partial_\sigma\tilde\eta=((\nabla u_\infty)\circ\tilde\eta)\partial_\sigma\tilde\eta+(\partial_\sigma U)\circ\tilde\eta
+(\nabla U\circ\tilde\eta)\partial_\sigma\tilde\eta
\end{equation*}
and then
\begin{equation}\label{key-equation}
\partial_t\partial_\sigma\tilde\eta|_{\sigma=0}=((\nabla u_\infty)\circ\tilde\eta)\partial_\sigma\tilde\eta|_{\sigma=0}+(\partial_\sigma U)\circ\tilde\eta|_{\sigma=0}.
\end{equation}
We see that 
\begin{equation*}
((\nabla u_\infty)\circ\tilde\eta)\partial_\sigma\tilde\eta=
\begin{pmatrix}
\partial_\sigma\tilde\eta_2(\partial_{2}f)\circ\tilde \eta_2\\
0
\end{pmatrix},
\end{equation*}
thus multiplying $\partial_\sigma\tilde \eta_2$ to  the second component of \eqref{key-equation}, and 
take integration in space, we have 
\begin{equation*}
(1/2)\partial_t\|\partial_\sigma\tilde\eta_2(0,t)\|_{L^2}^2
\leq \|\partial_\sigma U_2(0,t)\|_{L^2}\|\partial_\sigma\tilde\eta_2(0,t)\|_{L^2},
\end{equation*}
and this means
\begin{equation*}
\|\partial_\sigma\tilde\eta_2(0,t)\|_{L^2}\lesssim \int_0^t\|\partial_\sigma U_2(0,t')\|_{L^2}dt'.
\end{equation*}
We plug the above estimate into the first component of \eqref{key-equation}, we obtain 
\begin{equation*}
\begin{split}
(1/2)\partial_t\|\partial_\sigma\tilde\eta_1(0,t)\|_{L^2}^2\leq 
&\int^t_0\|\partial_\sigma U_2(0,t')\|_{L^2}dt'\|f'_m\|_{L^\infty}\|\partial_\sigma\tilde\eta_1(0,t)\|_{L^2}\\
+
&\|\partial_\sigma U_1(0,t)\|_{L^2}\|\partial_\sigma\tilde\eta_1(0,t)\|_{L^2},
\end{split}
\end{equation*}
and this means
\begin{equation*}
\|\partial_\sigma\tilde\eta_1(0,t)\|_{L^2}\lesssim \int^{t}_0\int^{t'}_0\|\partial_\sigma U_2(0,t'')\|_{L^2}dt''dt'\|f'_m\|_{L^\infty}+t\sup_{0<t'<t}\|\partial_\sigma U_1(0,t')\|_{L^2}.
\end{equation*}
On the other hand, 
 observe that by construction $\partial_\sigma\tilde{\eta}(\sigma, t)|_{\sigma=0}$ is a Jacobi field along the geodesic $\tilde{\eta}(0,t)$ and therefore we have 
(see \cite{M}, and, for the recent development on the inviscid parallel shear flow on compact manifolds, see \cite{TY})
\begin{equation*}
\|\partial_\sigma\tilde\eta(0,t)\|_{L^2}\geq Ct\quad \text{for any}\quad t>0,
\end{equation*}
 where  $C=\partial_t|_{t=0}\|\partial_\sigma\tilde\eta(0,t)\|_{L^2}$.
In our case this positive constant $C$ is nothing more than $\|U_{in}\|_{L^2}$, which is (essentially) independent of $m$.
In fact,
since
\begin{equation*}
\frac{\partial_\sigma\tilde\eta(0,t)}{t}\sim\frac{1}{t}\int_0^t (\partial_\sigma U\circ\tilde\eta)(0,t')dt'\sim U_{in},
\end{equation*}
from \eqref{key-equation}, we see that, at least $\partial_t|_{t=0}\|\partial_\sigma\tilde\eta(0,t)\|_{L^2}$ is nonzero (for each $m$).
Second we show that $C=\|U_{in}\|_{L^2}$ rigorously.
By
\begin{equation*}
\partial_t\|\partial_\sigma\tilde\eta(0,t)\|_{L^2}=\frac{\langle\partial_t\partial_\sigma\tilde\eta(0,t),\partial_\sigma\tilde\eta(0,t)\rangle_{L^2}}{\|\partial_\sigma\tilde\eta(0,t)\|_{L^2}}
\end{equation*}
and since the Jacobi field vanishes at $t=0$, i.e., $\partial_\sigma\tilde{\eta}(0,0)=0$, then, by L'H\^opital's rule, we have 
\begin{equation*}
\partial_t|_{t=0}\|\partial_\sigma\tilde\eta(0,t)\|_{L^2}=\frac{\langle\partial_t|_{t=0}\partial_\sigma\tilde\eta(0,t),\partial_t|_{t=0}\partial_\sigma\tilde\eta(0,t)\rangle_{L^2}}{\partial_t|_{t=0}\|\partial_\sigma\tilde\eta(0,t)\|_{L^2}}=\frac{\|U_{in}\|_{L^2}^2}{\partial_t|_{t=0}\|\partial_\sigma\tilde\eta(0,t)\|_{L^2}}
\end{equation*}
Thus we have $\partial_t|_{t=0}\|\partial_\sigma\tilde\eta(0,t)\|_{L^2}=\|U_{in}\|_{L^2}$.

Therefore we have 
\begin{equation*}
\begin{split}
 Ct\leq &\|\partial_\sigma\tilde\eta(0,t)\|_{L^2}\sim \|\partial_\sigma\tilde\eta_1(0,t)\|_{L^2}+\|\partial_\sigma\tilde\eta_2(0,t)\|_{L^2}\\
\lesssim &
\int_0^t\|\partial_\sigma U_2(0,t')\|_{L^2}dt'+\int^{t}_0\int^{t'}_0\|\partial_\sigma U_2(0,t'')\|_{L^2}dt''dt'\|f'_m\|_{L^\infty}\\
+&
t\sup_{0<t'<t}\|\partial_\sigma U_1(0,t')\|_{L^2}.
\end{split}
\end{equation*}
For any sufficiently small $\epsilon>0$, let us assume that both $\partial_\sigma U_1$ and $\partial_\sigma U_2$ are controlled by
\begin{equation*}
\|\partial_\sigma U_1(0,t)\|_{L^2}\lesssim 1 \quad\text{and}\quad \|\partial_\sigma U_2(\sigma,t)\|_{L^2}
\lesssim \langle t\rangle^{-1-\epsilon}
\end{equation*}
for $t>0$. Then this assumption breaks the lower bound estimate, and 
it is in contradiction.
Finally, for any sufficiently large $T>0$,
there is $\sigma_0>0$ such that for any $\sigma<\sigma_0$ and  $t\in[0,T]$,
\begin{equation*}
\|U_1(\sigma,t)\|_{L^2}\sim \sigma \|U_1(0,t)\|_{L^2}\gtrsim \sigma \quad\text{or}\quad \|U_2(\sigma,t)\|_{L^2}\sim \sigma\|U_2(0,t)\|_{L^2}\gtrsim \sigma \langle t\rangle^{-1-\epsilon}
\end{equation*}
by the 
Taylor expansion along $\sigma$-direction.


\vspace{0.5cm}
\noindent
{\bf Acknowledgments.}\ 
The author would like to thank Professors Herbert Koch and Gerard Misio{\l}ek for inspiring communications.
Research of TY  was partially supported by 
Grant-in-Aid for Young Scientists A (17H04825),
Grant-in-Aid for Scientific Research B (15H03621, 17H02860, 18H01136 and 18H01135).

\bibliographystyle{amsplain}

\end{document}